\documentclass{amsart}

\usepackage{amscd}

\begin{document}

\newtheorem{thm}{Theorem}[section]
\newtheorem{lem}[thm]{Lemma}
\newtheorem{remark}[thm]{Remark}
\newtheorem{cor}[thm]{Corollary}
\newtheorem{prop}[thm]{Proposition}
\theoremstyle{definition}
\newtheorem{ex}[thm]{Example}
\newtheorem*{question}{Question}
\newtheorem{defn}[thm]{Definition}
\newtheorem*{rem}{Remark}
\newtheorem*{ack}{Acknowledgements}
\newtheorem{eg}[thm]{Example}
\newtheorem{conj}[thm]{Conjecture}

\newcommand{\Q}{{\mathbb Q}}
\newcommand{\C}{{\mathbb C}}
\newcommand{\R}{{\mathbb R}}
\newcommand{\T}{{\mathbb T}}
\newcommand{\Tate}{{\rm Tate}}
\newcommand{\Th}{{\rm Th}}
\newcommand{\Z}{{\mathbb Z}}
\newcommand{\hT}{{\hat{\mathbb T}}}
\newcommand{\bT}{{\bf T}}
\newcommand{\cL}{\mathbb{L}}
\newcommand{\cH}{\mathbb{H}}
\newcommand{\gU}{\mathcal{U}}
\newcommand{\gG}{\mathcal{G}}

\title {Thom prospectra for loopgroup representations}

\author{Nitu Kitchloo and Jack Morava}

\address{Department of Mathematics, Johns Hopkins University,
Baltimore, Maryland 21218}
\email{nitu@math.jhu.edu, jack@math.jhu.edu}

\date {1 May 2004}

\begin{abstract} This is very much an account of work in progress.
We sketch the construction of an Atiyah dual (in the category of 
$\T$-spaces) for the free loopspace of a manifold; the main technical 
tool is a kind of Tits building for loop groups, discussed in detail 
in an appendix. Together with a new localization theorem for 
$\T$-equivariant $K$-theory, this yields a construction of the 
elliptic genus in the string topology framework of Chas-Sullivan, 
Cohen-Jones, Dwyer, Klein, and others. We also show how the Tits building 
can be used to construct the dualizing spectrum of the loop group. Using
a tentative notion of equivariant $K$-theory for loop groups, we relate 
the equivariant $K$-theory of the dualizing spectrum to recent work of 
Freed, Hopkins and Teleman. 
\end{abstract}

\maketitle

\begin{center}
{\bf Introduction}
\end{center} 

\noindent
If $P \to M$ is a principal bundle with structure group $G$ then 
$LP \to LM$ is a principal bundle with structure group 
\[
LG = {\rm Maps}(S^1,G) , 
\]
and if the tangent bundle of $M$ is defined by a representation $V$ of 
$G$ then the tangent bundle of $LM$ is defined by the representation
$LV$ of $LG$. The circle group $\T$ acts on all these spaces. \bigskip

\noindent
This is a report on the beginnings of a theory of differential topology 
for such objects. Note that if we want the structure group $LG$ to be 
connected, we need $G$ to be 1-connected; thus SU($n$) is preferable to
U($n$). This helps explain why Calabi-Yau manifolds are so central in 
string theory, and this note is written assuming this simplifying 
hypothesis. \bigskip

\noindent
Alternately, we could work over the universal cover of $LM$; then
$\pi_2(M)$ would act on everything by decktranslations, and our 
topological invariants become modules over the Novikov ring $\Z[H_2(M)]$.
From the point of view we're developing, these translations seem to be 
what really underlies modularity, but this issue, like several others,
will be backgrounded here. \bigskip

\noindent
The circle action on the free loopspace defines a structure much closer
to classical differential geometry than one finds on more general (eg)
Hilbert manifolds; this action defines something like a Fourier filtration 
on the tangent space of this infinite-dimensional manifold, which is
in some sense locally finite. This leads to a host of new kinds of 
geometric invariants, such as the Witten genus; but this filtration
is unfamiliar, and has been difficult to work with \cite{CoS}. The main 
conceptual result of this note [which was motivated by ideas of Cohen, 
Godin, and Segal] is the definition of a canonical equivariant `thickening' 
of a free loopspace, where the pulled-back tangent bundle admits 
a canonical filtration by finite-dimensional equivariant bundles. This 
thickening involves a contractible $LG$-space called the affine 
Tits building ${\bf A}(LG)$. This space occurs under various guises in nature:
it is a homotopy colimit of homogeneous spaces with respect to a finite 
collection of compact Lie subgroups of $LG$. It is also the affine space 
of principal $G$-connections on the trivial bundle over $S^1$. We explore 
its structure in the appendix. \bigskip

\noindent
In the first section below, we recall why the Spanier-Whitehead dual
of a finite CW-space is a ring-spectrum, and sketch the construction
(due to Milnor and Spanier, and Atiyah) of a model for that dual, when
the space is a smooth compact manifold. Our goal is to produce an 
analog of this construction for a free loopspace, which captures as
much as possible of its string-topological algebraic structure. In the
second section, we introduce the technology used in our construction:
pro-spectra associated to filtered infinite-dimensional vector bundles,
and the topological Tits building which leads to the construction
of such a filtration for the tangent bundle. \bigskip

\noindent
In \S 3 we observe that recent work of Freed, Hopkins, and Teleman on the 
Verlinde algebra can be reformulated as a conjectural duality between
$LG$-equivariant K-theory of a certain dualizing spectrum
for $LG$ constructed from its Tits building, and positive-energy 
representations of $LG$. In \S 4 we use a new strong localization 
theorem to study the equivariant $K$-theory of our construction, and we show 
how this recovers the Witten genus from a string-topological point of view. 
\bigskip

\noindent
We plan to discuss actions of various string-topological operads \cite{CJ} 
on our construction in a later paper; that work is in progress. \bigskip

\noindent
We would like to thank R.\ Cohen, V.\ Godin, and A.\ Stacey for 
many helpful conversations, and we would also like to acknowledge 
work of G.\ Segal and S.\ Mitchell as motivation for many of the 
ideas in this paper. \bigskip

\section{The Atiyah dual of a manifold}

\noindent
If $X$ is a finite complex, then the function spectrum 
$F(X,S^0)$ is a {\bf ring-spectrum} (because $S^0$ is).
If $X$ is a manifold $M$, Spanier-Whitehead duality says that 
\[
F(M_+,S^0) \sim M^{-TM} \;.
\]
If $E \to X$ is a vector bundle over a compact space, we can define its
Thom space to be the one-point compactification
\[
X^E := E_+ \;.
\]
There is always a vector bundle $E_\perp$ over $X$ such that 
\[
E \oplus E_\perp \cong {\bf 1}_N
\]
is trivial, and following Atiyah, we write
\[
X^{-E} := S^{-N}X^{E_\perp} .
\]
With this notation, the Thom collapse map for an embedding $M \subset \R^N$
is a map
\[
S^N = \R^N_+ \to M^\nu = S^N \; M^{-TM} ,
\]
and the midpoint construction
\[
M_+ \wedge M^\nu \to \R^N_+ = S^N
\]
defines the equivalence with the functional dual. More generally, a smooth map
$f: M \to N$ of compact closed orientable manifolds has a Pontrjagin-Thom 
dual map
\[
f_{PT} : N^{-TN} \to M^{-TM} 
\]
of spectra; in particular, the map $S^0 \to M^{-TM}$ dual to the projection 
to a point defines a kind of fundamental class, and the dual to the diagonal 
of $M$ makes $M^{-TM}$ into a ring-spectrum \cite{C}. \bigskip

\noindent
{\bf Prospectus:} Chas and Sullivan \cite{CS} have constructed a very 
interesting product on the homology of a free loopspace, suitably desuspended,
motivated by string theory. Cohen and Jones \cite{CJ} saw that this product 
comes from a ring-spectrum structure on 
\[
LM^{-TM} := LM^{-e^*TM}
\]
where 
\[
e : LM \to M
\]
is the evaluation map at $1 \in S^1$. Unfortunately this evaluation map is 
not $\T$-equivariant, so the Chas-Sullivan Cohen-Jones spectrum is not in 
general a $\T$-spectrum. The full Atiyah dual constructed below promises to 
capture some of this equivariant structure. The Chas-Sullivan Cohen-Johes 
spectrum and the full Atiyah dual live in rather different worlds: our 
prospectrum is an equivariant object, whose multiplicative properties are 
not year clear, while the CSCJ spectrum has good multiplicative properties, 
but it is not a $\T$-spectrum. In some vague sense our object resembles 
a kind of center for the Chas-Sullivan-Cohen-Jones spectrum, and we hope that
a better understanding of the relation between open and closed strings
will make it possible to say something more explict about this. 
\bigskip

\section{Problems \& Solutions}

\noindent
For our constructions, we need two pieces of technology:\bigskip

\noindent
Cohen, Jones, and Segal \cite{CJS}(appendix) associate to a filtration
\[
{\bf E} : \dots \subset E_i \subset E_{i+1} \subset \dots
\]
of an infinite-dimensional vector bundle over $X$, a pro-object
\[
X^{-{\bf E}} : \dots \to X^{-E_{i+1}} \to X^{-E_i} \to \dots
\]
in the category of spectra. [A rigid model for such an object can be
constructed by taking {\bf E} to be a bundle of Hilbert spaces, which 
are trivializable by Kuiper's theorem. Choose a trivialization ${\bf E} 
\cong H \times X$ and an exhaustive filtration $\{H_k\}$ of $H$ by
finite-dimensional vector spaces; then we can define
\[
X^{-E_i} = \lim S^{-H_k} X^{H_k \cap E_i^\perp} \;,
\]
with $E_i^\perp$ the orthogonal complement of $E_i$ in the trivialized
bundle {\bf E}.] This pro-object will, in general, depend on the choice 
of filtration. We will be interested in the {\bf direct} systems associated
to such a pro-object by a cohomology theory; of course in general the
colimit of this system can be very different from the cohomology of the
limit of the system of pro-objects.\bigskip

\begin{ex}
If $X = \C P_\infty$, $\eta$ is the Hopf bundle, and {\bf E} is
\[
\infty \eta : \dots \subset (k-1) \eta \subset k \eta \subset (k+1) \eta 
\subset \dots
\]
then the induced maps of cohomology groups are multiplication by the Euler 
classes of the bundles $E_{i+1}/E_i$, so 
\[
H^*(\C P_\infty^{-\infty \eta},\Z) := {\rm colim} \{\Z[t],t-{\rm mult}\} = 
\Z[t,t^{-1}] \;.
\] 
\end{ex}

\noindent
We would like to apply such a construction to the tangent bundle of
a free loopspace. Unfortunately, these tangent bundles do {\bf not}, in 
general, possess any such nice filtration by finite-dimensional 
($\T$-equivariant) subbundles \cite{CoS}! However, such a splitting
{\bf does} exist in a neighborhood of the {\bf constant} loops:
\[
M = LM^\T \subset LM
\]
has normal bundle 
\[
\nu (M \subset LM) = TM \otimes_\C (\C[q,q^{-1}]/\C)
\]
(at least, up to completions; and assuming things complex for convenience).
Here small perturbations of a constant loop are identified with their 
Fourier expansions
\[
\sum_{n \in \Z} a_n q^n ,
\]
with $q = e^{i \theta}$. The related fact, that $TLM$ is defined by the 
representation $LV$ of $LG$ looped up from the finite-dimensional 
representation $V$ of $G$, will be important below: for $LV$ is {\bf not} 
a positive-energy representation of $LG$. \bigskip

\noindent
The main step toward our resolution of this problem depends on
the following result, proved in \S 7 below. Such constructions were
first studied by Quillen, and were explored further by S. Mitchell 
\cite{Mi}. The first author has studied these buildings for 
a general Kac-Moody group \cite{Ki}; most of the properties of the 
affine building used below hold for this larger class. \bigskip

\begin{thm}
The {\bf topological affine Tits building} 
\[
{\bf A}(LG) := {\rm hocolim}_I  \; LG/H_I
\]
of $LG$ is $\T\tilde{\times} LG$-{\bf equivariantly} contractible. In other 
words, given any compact subgroup $K \subset \T\tilde{\times} LG$, the fixed 
point space ${\bf A}(LG)^K$ is contractible. \end{thm}

\noindent
[Here $I$ runs over certain proper subsets of roots of $G$, and the $H_I$ are
certain {\bf compact} `parabolic' subgroups of $LG$ (see \S 7.2).] 

\noindent
\begin{remark}
{\rm The group $LG$ admits a universal central extension $\cL G$. The 
natural action of the rotation group $\T$ on $LG$ lifts to $\cL G$, and
the $\T$-action preserves the subgroups $H_I$. Hence ${\bf A}(LG)$ admits an 
action of $\T \tilde{\times} \cL G$, with the center acting trivially. We 
can therefore express ${\bf A}(LG)$ as
\[ 
{\bf A}(LG) = {\rm hocolim}_I \; \cL G/\cH_I 
\]
where $\cH_I$ is the induced central extension of $H_I$.} \end{remark}

\section*{Other descriptions of ${\bf A}(LG)$}

\noindent
This Tits building has other descriptions as well. For example:

\medskip
\noindent {\bf 1.}
${\bf A}(LG)$ can be seen as the classifying space for proper 
actions with respect to the class of compact Lie subgroups of $\T\tilde{\times} 
\cL G$.

\medskip
\noindent {\bf 2.}
It also admits a more differential-geometric description as the smooth 
infinite dimensional manifold of holonomies on $S^1 \times G$ (see Appendix):
Let ${\mathcal S}$ denote the subset of the space of smooth maps from $\R$ 
to $G$ given by
\[ 
{\mathcal S} = \{ g(t):\R \rightarrow G, \, \, g(0)=1, \, \, g(t+1) = g(t)
\cdot g(1) \} 
\;;
\]
then ${\mathcal S}$ is homeomorphic to ${\bf A}(LG)$. The action of $h(t) 
\in LG$ on $g(t)$ is given by $hg(t) = h(t) \cdot g(t) \cdot h(0)^{-1}$, 
where we identify the circle with $\R/\Z$. The action of $ x \in \R/\Z = 
\T$ is given by $xg(t) = g(t+x) \cdot g(x)^{-1}$.  

\medskip 
\noindent {\bf 3.}
The description given above shows that ${\bf A}(LG)$ is equivalent to the 
affine space ${\mathcal A}(S^1\times G)$ of connections on the trivial 
$G$-bundle $S^1 \times G$. This identification associates to the function 
$f(t) \in {\mathcal S}$, the connection $f'(t)f(t)^{-1}$. Conversely, the 
connection $\nabla_t$ on $S^1 \times G$ defines the function $f(t)$ given by 
transporting the element $(0,1) \in \R \times G$ to the point $(t,f(t)) \in 
\R \times G$ using the connection $\nabla_t$ pulled back to the trivial 
bundle $\R \times G$. 

\begin{remark}
{\rm These equivalent descriptions have various useful consequences. 
For example, the model given by the space ${\mathcal S}$ of holonomies
says that given a finite cyclic group $H \subset \T$, the fixed point space 
${\mathcal S}^H$ is homeomorphic to ${\mathcal S}$. Moreover, this is a 
homeomorphism of $LG$-spaces, where we consider ${\mathcal S}^H$ as an 
$LG$-space and identify $LG$ with $LG^H$ in the obvious way. Notice also 
that ${\mathcal S}^{\T}$ is $G$-homeomorphic to the model of the adjoint
representation of $G$ defined by ${\rm Hom}(\R,G)$. \medskip

\noindent
Similarly, the map ${\mathcal S} \rightarrow G$ given by evaluation at $t=1$ 
is a principal $\Omega G$ bundle, and the action of $G = LG/\Omega G$ 
on the base $G$ is given by conjugation. This allows us to relate our work 
to that of Freed, Hopkins and Teleman in the following section. \medskip

\noindent
Finally, the description of ${\bf A}(LG)$ as the affine space 
${\mathcal A}(S^1 \times G)$ implies that the fixed point space 
${\bf A}(LG)^K$ is contractible for any compact subgroup $K \subseteq 
\T \tilde{\times} \cL G$.} 
\end{remark} \bigskip

\noindent
If $E \to B$ is a principal bundle with structure group $LG$, then 
(motivated by ideas of \cite{CG}) we construct a `thickening' of $B$: 

\begin{defn} {\it The thickening of $B$ associated to the bundle $E$ is
the $LG$-space}
\[ 
B_!(E) = E \times_{LG} {\bf A}(LG) = {\rm hocolim}_I \; E/H_I \;.
\]
We will omit $E$ from the notation, when the defining bundle is clear
from context. \end{defn}

\begin{remark}
{\rm If $P \to M$ is a principal $G$ bundle, then $LP \to LM$ is a principal
$LG$ bundle. In this case, the description above gives $L_!M := LM_!(LP)$ 
a smooth structure:
\[ 
L_!M = \{ (\gamma,\omega) \, | \quad \gamma \in LM, \quad \omega \in 
{\mathcal A}(\gamma^*(P)) \} 
\]
where ${\mathcal A}(\gamma^*(P))$ is the space of connections on the 
pullback bundle $\gamma^*(P)$.}
\end{remark} \bigskip

\noindent
Let $\T \tilde{\times} \cL G$ be the extension of the central extension of 
$LG$ by $\T$ acting as rotations. On restriction to the subgroup 
$\T \tilde{\times} \cH_I$, a unitary representation $U$ of $\T \tilde{\times} 
\cL G$ decomposes into a sum of finite dimensional representations. We want 
to construct a Thom $\T \tilde{\times} \cL G$-prospectrum ${\bf A}
(LG)^{-U}$. \bigskip

\noindent  
We consider the decomposition of the restriction of $U$ to $\T \tilde{\times} 
\cH_I$ as a sum of irreducibles:
\[
U|_{\T \tilde{\times} \cH_I} \cong \oplus \; U_I(\alpha)
\]
and let 
\[
U_I(k) = \oplus \; \{ U_I(\alpha) \;|\; \dim U_I(\alpha) \leq k\} \;.
\]
Then we can define $S_I^{-U}$ to be the Thom $\T \tilde{\times} 
\cH_I$-prospectrum associated to the filtered (equivariant) vector bundle
\[
{\bf U}_I : \dots \subset U_I(k) \subset U_I(k+1) \subset \dots 
\]
over a point. If $I \subset J$ then $\cH_I$ maps naturally to $\cH_J$, and 
there is a corresponding morphism
\[
{\bf U}_J \to {\bf U}_I
\]
of filtered vector bundles, given by inclusions $U_J(k) \to U_I(k)$. 
\bigskip

\begin{defn} {\it We define ${\bf A}(LG)^{-U}$ to be the $\T \tilde{\times} 
\cL G$-prospectrum
\[ 
{\bf A}(LG)^{-U} = {\rm hocolim}_I \; \cL G_+ \wedge_{\cH_I} S_I^{-U} \;,
\]
where $\cL G_+$ denotes $\cL G$, with a disjoint basepoint.} 
\end{defn}

\noindent
Homotopy colimits in the category of prospectra can be defined in 
general, using the model category structure of \cite{CI}. \bigskip

\begin{remark}
{\rm Given any principal $LG$-bundle $E \rightarrow B$, and a
representation $U$ of $LG$, we define the Thom prospectrum of the 
virtual bundle associated to the representation $-U$ to be}
\[ 
B_!^{-U} = E_+ \wedge_{LG} {\bf A}(LG)^{-U} = {\rm hocolim}_I
\; E_+\wedge_{H_I}S_I^{-U} \;.  
\]
\end{remark}

\noindent
In particular, if $P$ is the refinement of the frame bundle of $M$ via a 
representation $V$ of $G$, then the tangent bundle of $LM$ is defined by 
the representation $LV$ of $LG$. 

\begin{defn}
{\it The Atiyah dual $LM^{-{\bf T}LM}$ of $LM$ is the pro-spectrum 
$L_!M^{-LV}$.} \end{defn}

\noindent
We will explore this object further in \S 6.

\section{The dualizing spectrum of $\cL G$}

\noindent
The dualizing spectrum of a topological group $K$ is defined \cite{Klein} 
as the $K$-homotopy fixed point spectrum: 
\[ 
D_K = K_+^{hK} = F(EK_+,K_+)^K  
\]
where $K_+$ is the suspension spectrum of the space $K_+$, endowed with a 
right $K$-action. The dualizing spectrum $D_K$ admits a $K$-action given by 
the residual left $K$-action on $K_+$. If $K$ is a compact Lie group, then 
it is known that $D_K$ is the one point compactification of the adjoint 
representation $S^{Ad(K)}$. It is also known that there is a $K \times 
K$-equivariant homotopy equivalence 
\[ 
K_+ \cong F(K_+,D_K) \;.
\]
It follows from the compactness of $K_+$ that for any free $K_+$-spectrum 
$E$, we have the $K$-equivariant homotopy equivalence 
\[ 
E \cong F(K_+,E\wedge_{K_+} D_K) \;.
\]
It is our plan to understand the dualizing spectrum for the (central 
extension of the) loop group. 

\begin{thm} There is an equivalence 
\[ 
D_{LG} \cong {\rm holim}_I \; LG_+\wedge_{H_I} S^{Ad(H_I)} 
\]
of left $LG$-spectra.
\end{thm}

\begin{proof}
We have the sequence of equivalences:
\[ 
D_{LG} = F(ELG_+,LG_+)^{LG} \cong F(ELG_+ \wedge {\bf A}(LG)_+,LG_+)^{LG} \;.
\]
The final space may be written as
\[ 
{\rm holim}_I \; F(ELG_+ \wedge_{H_I} LG_+, LG_+)^{LG} = {\rm holim}_I \; 
F(ELG_+, LG_+)^{H_I} \;.
\]
Now recall the equivalence of $H_I \times H_I$-spectra:
\begin{equation} \label{equiv}
LG_+ \cong F({H_I}_+,LG_+\wedge_{H_I} D_{H_I}) \;.
\end{equation}
Taking $H_I$-homotopy fixed points implies a left $H_I$-equivalence
\[ 
F(ELG_+, LG_+)^{H_I} = (LG_+)^{hH_I} \cong LG_+ \wedge_{H_I} S^{Ad(H_I)} \;;
\]
where we have used equation (\ref{equiv}) at the end. Replacing 
this term into the homotopy limit completes the proof. \end{proof}

\noindent
Similarly, we have:
\begin{thm} There is an equivalence
\[ 
D_{\cL G} \cong {\rm holim}_I \; \cL G_+\wedge_{\cH_I} S^{Ad(\cH_I)} 
\]
of left $\cL G$-spectra.
\end{thm}

\begin{remark}
{\rm The diagram underlying $D_{LG}$ or $D_{\cL G}$ can be constructed in the
category of spaces. Given an inclusion $I \subseteq J$, the orbit of a 
suitable element in $Ad(\cH_J)$ gives an embedding $\cH_J/\cH_I \subset 
Ad(\cH_J)$, and the Pontrjagin-Thom construction for this embedding defines
an $\cH_J$-equivariant map 
\[ 
S^{Ad(\cH_J)} \longrightarrow {\cH_J}_+\wedge_{\cH_I}S^{Ad(\cH_I)} 
\]
which extends to the map
\[  
\cL G_+\wedge_{\cH_J} S^{Ad(\cH_J)} \longrightarrow \cL G_+\wedge_{\cH_I} 
S^{Ad(\cH_I)}
\]
required for the diagram. Moreover, composites of these maps can be made 
compatible up to homotopy.}
\end{remark} 

\section*{A conjectural relationship with the work of Freed, Hopkins and 
Teleman}
 
\noindent
The discussion below assumes the existence of a hypothetical
${\cL G}$-equivariant $K$-theory, whose value on a point in degree zero is 
the Grothendieck group of positive energy representations (or equivalently,
the group of characters of integrable representations). The symmetric 
monoidal category whose objects are finite direct sums of irreducible 
positive energy representations, and whose morphism spaces consist of the 
(nonequivariant) isomorphisms of the vector spaces underlying the
representation (given the compactly generated topology) defines a candidate 
for a spectrum representing such a functor: this is a topological category 
with an ${\cL G}$-action which respects the symmetric monoidal structure. 
\bigskip

\noindent
This hypothesis provides us with a convenient language. We expect to return
to the underlying technical issues in a later paper.\bigskip

\noindent
The center of $\cL G$ acts trivially on $D_{\cL G}$, defining a second
grading on $K_{\cL G}^*(D_{\cL G})$; we will use a formal variable $z$ to 
keep track of the grading, so
\[ 
K_{\cL G}^*(D_{\cL G}) = \bigoplus_n K_{\cL G}^{*,n}(D_{\cL G}) z^n. 
\]
The spectral sequence for the cohomology 
of a cosimplicial space, in the case of $K_{\cL G}^*(D_{\cL G})$, has
\[ 
E_2^{i,j} = {\rm colim}^i_I \; K^j_{\cH_I}(S^{Ad(\cH_I)}) \;.
\]
This spectral sequence respects the second grading given by powers of $z$. 
In a sequel to this paper, we will show that this spectral sequence collapses 
to give 
\[ 
K_{\cL G}^*(D_{\cL G}) = {\rm colim}_I \; K^*_{\cH_I}(S^{Ad(\cH_I)}) \cong 
{\rm colim}_I \; K^{*-r-1}_{\cH_I}(pt) \;, 
\]
where $r$ is the rank of $G$. Therefore, this group admits a natural Thom 
class given the system $\{ S(Ad(\cH_I))\}$ of spinor bundles for the adjoint 
representations of the parabolics $\cH_I$. In section 11 of \cite{FHT}, the 
authors construct an explicit map between the Verlinde algebra and this colimit, as follows: \medskip

\noindent
To a positive energy representation corresponding to a dominant character 
$\lambda$, we associate the $\cL G$-equivariant bundle given by ${\mathcal 
L}_{-\lambda -\rho}\otimes S(N)$, where ${\mathcal L}_{-\lambda -\rho}$ is 
the canonical line bundle over the coadjoint orbit of the regular element 
$\lambda + \rho$, and $S(N)$ is the spinor bundle of the normal bundle to 
the coadjoint orbit. Such an orbit is of the form $\cL G/\cH$ for some 
parabolic subgroup $\cH$, and its normal bundle is $Ad(\cH)$, so this element 
defines a class in $K_{\cL G}(D_{\cL G})$. The same can be done for 
antidominant weights. This suggests the following:

\begin{conj}
For following map is an isomorphism in homogeneous degree $z^n$, for $n\neq 
0$: 
\[ 
\bigoplus_{k\geq 0} V_{k}z^{k+h} \; \bigoplus_{k\geq 0} V_{k}z^{-(k+h)} \cong 
{\rm colim}_I \; K_{{\cH}_I}(pt) \rightarrow K_{\cL G}^{r+1}(D_{\cL G}) = 
\bigoplus_n K_{\cL G}^{*,n}(D_{\cL G}) z^n
\]
where $V_k$ is the Verlinde algebra of level $k$, $h$ is the dual Coxeter 
number of $G$, and $r$ is its rank. 
\end{conj}

\begin{eg} To illustrate this in an example, consider the case $G = SU(2)$. 
In this case $r=1$, $h(G)=2$. Here the groups 
${\cH}_I$ are given by 
\[ 
{\cH}_0 = SU(2) \times S^1, \quad {\cH}_1 = S^1 \times SU(2), \quad {\cH}_0 
\cap {\cH}_1 = T = S^1 \times S^1 \;.
\]
The respective representation rings may be identified by restriction with 
subalgebras of $K_T(pt) = \Z[u^{\pm1},z^{\pm1}]$ :
\[ 
K_{{\cH}_0}(pt) = \Z[u+u^{-1},(z/u)^{\pm1}], \quad K_{{\cH}_1}(pt) = 
\Z[z^{\pm1}, u+u^{-1}] \;.
\]
Now consider the two pushforward maps involved in the colimit:
\[ 
\varphi_0:K_T(pt) \rightarrow K_{{\cH}_0}(pt), \quad \varphi_1:K_T(pt) 
\rightarrow K_{{\cH}_1}(pt) 
\]
A quick calculation shows that for $k > 0$, we have 
\[ 
\varphi_j(z^k) = \begin{cases} (z/u)^{k}Sym^k(u+u^{-1}), \quad j=0 \\
z^k , \quad \quad \quad \quad \quad \quad \quad \quad \quad \quad j=1  
\end{cases} 
\]
\[ \varphi_j(z^ku^{-1}) = \begin{cases} (z/u)^kSym^{k-1}(u+u^{-1}), 
\quad j=0 \\ 
0, \quad \quad \quad \quad \quad \quad \quad \quad \quad \quad \quad  j=1 \;, 
\end{cases} 
\]
where $Sym^k(V)$ denotes the $k$-th symmetric power of the representation $V$,
e.g. $Sym^k(u + u^{-1}) = u^k + \cdots + u^{-k}$. \bigskip

\noindent
We also have a similar formula for negative exponents:
\[ 
\varphi_j(z^{-k}) = \begin{cases} -(u/z)^kSym^{k-2}(u+u^{-1}), \quad j=0 \\
z^{-k} ,  \quad \quad \quad \quad \quad \quad \quad \quad \quad \quad \quad 
j=1  \end{cases} 
\]
\[ \varphi_j(z^{-k}u^{-1}) = \begin{cases} -(u/z)^kSym^{k-1}(u+u^{-1}), 
\quad j=0 \\ 
0, \quad \quad \quad \quad \quad \quad \quad \quad \quad \quad \quad \quad  
j=1 \;. \end{cases}
\]
The colimit is the cokernel of
\[ 
\varphi_1\oplus \varphi_0 : K_T(pt) \longrightarrow K_{\cH_1}(pt) \oplus 
K_{\cH_0}(pt) 
\]
Now consider the decomposition 
\[ 
\Z[u^{\pm1},z^{\pm1}] = \Z[u+u^{-1},z^{\pm1}] \oplus u^{-1}\Z[u+u^{-1},
z^{\pm1}] \;.
\]
It is easy to check from this that the cokernel for nontrivial powers of $z$
is isomorphic to the cokernel of $\varphi_0$ restricted to $u^{-1}\Z[u+u^{-1},
z^{\pm1}]$ and hence is
\[ 
\bigoplus_{k\geq 0} \frac{\Z[u+u^{-1}]}{\langle Sym^{k+1}(u+u^{-1})\rangle} 
(z/u)^{k+2} \bigoplus_{k\geq 0} \frac{\Z[u+u^{-1}]}{\langle Sym^{k+1}
(u+u^{-1})\rangle} (u/z)^{k+2}
\] 
which agrees with the classical result \cite{Fr}. 
\end{eg}

\begin{remark}
{\rm We can calculate the equivariant $K$-homology ${K_{\cL G}}_*({\bf A}
(LG))$ using the same spectral sequence. This establishes an isomorphism 
between $K^*_{\cL G}(D_{\cL G})$ and ${K_{\cL G}}_*({\bf A}(LG))$. Results of
\cite{FHT} suggest that the latter group calculates the Verlinde algebra, 
which is yet another motivation for the conjecture. Recall also that 
${\bf A}(LG)$ is the classifying space for {\bf proper} actions (i.e. with 
compact isotropy) so our conjecture is a topological analog of the 
Baum-Connes conjecture for finite groups \cite{BK}} \end{remark}

\begin{question}
Given a manifold $LM$, with frame bundle $LP$, we can construct a spectrum
\[ 
D_{LM} := {\rm holim}_I \; LP_+\wedge_{H_I} S^{Ad(H_I)} 
\]
It would be very interesting to understand something about $K_{\T}(D_{LM})$.
\end{question}

\section{Localization Theorems} 

\noindent
If $E$ is a $\T$-equivariant complex-oriented multiplicative
cohomology theory, and $X$ is a $\T$-space, we have contravariant 
($j^*$) and covariant ($j^!$) homomorphisms associated to the 
fixedpoint inclusion 
\[
j : X^\T \subset X \;,
\]
satisfying
\[
j^* j^!(x) = x \cdot e_\T(\nu) \;;
\]
if the Euler class of the normal bundle $\nu$ is invertible, this leads to
a close relation between the cohomology of $X$ and $X^\T$. \bigskip

\noindent
More generally, if $f : M \to N$ is an equivariant map, then its 
Pontrjagin-Thom transfer is related to the analogous transfer defined 
by its restriction
\[
f^\T : M^\T \to N^\T
\]
to the fixedpoint spaces, by a `clean intersection' formula:
\[
j^*_N \circ f^! (-) = f^{\T!}(j^*_M (-) \cdot e_\T(\nu(f)|_{M^\T})) \;.
\]

\begin{defn}
{\it The fixed-point orientation defined by the Thom class
\[
\Th^\dagger(\nu(f^\T)) = \Th(\nu(f^\T)) \cdot e_\T(\nu(f)|_{M^\T}
\]
for the normal bundle of the inclusion of fixed-point spaces is
the product of the usual Thom class with the equivariant Euler 
class of the full normal bundle restricted to the fixed-point 
space.} \end{defn}

\noindent
Since $f^{\T!}(-) = f_{PT}^{\T*}(- \cdot \Th(\nu(f^\T)))$, in this
new notation the clean intersection formula becomes 
\[
j^*_N \circ f^! = f^{\T\dagger} \circ j^*_M 
\]
with a new Pontrjagin-Thom transfer
\[
f^{\T\dagger}(-) = f_{PT}^{\T*}( - \cdot \Th^\dagger(\nu(f^\T))) \;.
\]
In the case of most interest to us (free loopspaces), we identified
the normal bundle above, in \S 2; using that description, we have
\[
e_\T(\nu(M \subset LM)) = \prod_{0 \neq k \in \Z;i} (e(L_i)+_E [k](q)) ,
\]
where the $L_i$ are the line bundles in a formal decomposition of $TM$, 
$q$ is the Euler class of the standard one-dimensional complex 
representation of $\T$, and $+_E$ is the sum with respect to the formal 
group law defined by the orientation of $E$. It may not be immediately
obvious, but it turns out that such a formula implies that
the fixed-point orientation defined above will have good multiplicative
properties. \bigskip

\noindent
Such Weierstrass products sometimes behave better when `renormalized',
by dividing by their values on constant bundles \cite{AM}. If $E$ is 
$K_\T$ with the usual complex (Todd) orientation, we have
\[
e(L)+_K [k](q) = 1 - q^k L \;;
\]
but for our purposes things turn out better with the Atiyah-Bott-Shapiro 
{\bf spin} orientation; in that case the corresponding Euler class is
\[
(q^kL)^{1/2} - (q^kL)^{-1/2} \;.
\]
The square roots make sense under the simple-connectivity
assumptions on $G$ mentioned in the introduction: \bigskip

\noindent
To be precise, let $V$ be a representation of $LG$ with an intertwining 
action of $\T$. We restrict ourselves to representations $V$ which (for 
lack of a better name {\bf [?]}) we call {\bf symmetric}, i.e. such 
that $V$ is equivalent to the representation of $LG$ obtained by composing 
$V$ with the involution of $LG$ which reverses the orientation of the loops. 
The restriction of the representation $V$ to the constant loops $\T \times G 
\subset \T \tilde{\times} LG$ has a decomposition
\[ 
V = V^{\T} \oplus \sum_{k\neq 0} \; V_k \; q^k 
\]
where $V_k$ are representations of $G$, and $q$ denotes the fundamental 
representation of $\T$. Let $V(m)$ be the finite dimensional subrepresentation
\[ 
V(m) = V^{\T} \oplus \sum_{0 < |k| \leq m}V_k \; q^k \;;
\]
the symmetry assumption implies that $V_k = V_{-k}$ as representations of 
$G$, so this can be rewritten
\[
V(m) = V^{\T} \oplus \sum_{0 < k \leq m}V_k \; (q^k\oplus q^{-k}) \;.
\]
At this point we need the following \bigskip

\begin{prop} If $G$ is a compact Lie group, and $W$ is an 
$m$-dimensional complex spin representation of $G$, then the representation 
$\tilde{W} = W \otimes (q^k\oplus q^{-k})$ of $\T \times G$ admits a 
canonical spin structure. 
\end{prop}

\begin{proof} The representation $q^k \oplus q^{-k}$ of $\T$ admits a 
unique spin structure. Since $W$ is also endowed with a spin 
structure, the representation $W \otimes (q^k \oplus q^{-k})$ admits a 
canonical spin structure defined by their tensor product. 
\end{proof}

\begin{remark} If $G$ is simply connected, then any representation $W$ of 
$G$ admits a unique spin structure. 
\end{remark}

\noindent
This justifies the square roots of the formal line bundles appearing in 
the restriction of $TLM$ to $M$. The resulting renormalized Euler class
\[
\prod_{k \neq 0} \frac{(q^kL)^{1/2} - (q^kL)^{-1/2}}{q^{k/2} - q^{-k/2}}
\]
is a product of terms of the form
\[
(1-q^kL)(-q^kL)^{-1/2}(q^{-k}L)^{1/2}(1-q^kL^{-1})
\]
divided by terms of the form
\[
(1-q^k)(-q^k)^{-1/2}(q^{-k})^{1/2}(1-q^k) \;,
\]
(where now all $k$'s are positive) yielding a unit
\[
\epsilon_\T(L) = \prod_{k \geq 1} \frac{(1 - q^k L)(1 - q^k L^{-1})}
{(1 - q^k)^2} 
\]
in the ring $\Z[L^\pm][[q]]$.\bigskip

\noindent
Following 4.1, we can reformulate the localization theorem in 
terms of a new orientation, obtained by multiplying the ABS Thom class
by this unit, to get {\bf precisely} the Mazur-Tate normalization 
\[
\sigma(L,q) = (L^{1/2} - L^{-1/2})\prod_{k \geq 1} \frac{(1 - q^k L)(1 - q^k 
L^{-1})}{(1 - q^k)^2}
\]
for the Weierstrass sigma-function as Thom class for a line bundle $L$. 
This extends by the splitting principle to define the orientation giving
the Witten genus \cite{Za}. 
\bigskip

\section{One moral of the story}

\noindent
Since the early 80's physicists have been trying to interpret
\[
M \mapsto K_\T(LM)
\]
as a kind of elliptic cohomology theory; but of course we know better,
because we know that mapping-space constructions (such as free loop
spaces) don't preserve cofibrations. \bigskip

\noindent
Now it is an easy exercise in commutative algebra to prove that
\[
\Z((q)) := \Z[[q]][q^{-1}] 
\]
is flat over 
\[
K_\T = \Z[q^\pm] \;,
\]
for the completion of a Noetherian ring, eg $\Z[q]$, at an ideal, eg $(q)$, 
defines a flat \cite{AtM}(\S 10.14) $\Z[q]$-algebra $\Z[[q]]$. Flat modules 
pull back to flat modules \cite{B} (Ch I \S 2.7), so it follows that 
$\Z((q))$ is flat over the localization
\[
\Z[q^\pm] := \Z[q][q^{-1}] = \Z[q,q^{-1}] \;.
\]
The Weierstrass product above is a genuine formal power
series in $q$, so for questions involving the Witten genus it is
formally easier to work with the functor defined on finite $\T$-CW spaces by
\[
X \mapsto K^*_\T (X) \otimes_{\Z[q^{\pm}]} \Z((q)) := K^*_\hT (X) \;.
\]
This takes cofibrations to long exact sequence of $\Z((q))$-modules.
Its real virtue, however, is that it satisfies a strong localization 
theorem: \bigskip

\begin{thm}
If $X$ is a finite $\T$-CW space, then restriction to the 
fixedpoints defines an {\bf isomorphism}
\[
j^* : K_\hT (X) \cong K_\hT (X^\T) \;.
\]
\end{thm}

\noindent
Proof, by skeletal induction; based on the 
\begin{lem}
If $C \subset \T$ is a proper closed subgroup, then
\[
K^*_\hT (\T/C) = 0 \;.
\]
\end{lem}
\begin{proof}If $C$ is cyclic of order $n \neq 1$, then
\[
K_\hT (\T/C) = K_\T (\T/C) \otimes_{\Z[q^\pm]} \Z((q)) = \Z[q]/(q^n
-1) \otimes_{\Z[q]} \Z((q)) 
\]
is zero, since 
\[
- 1 = (q^n - 1) \cdot \sum_{k \geq 0} q^{nk} = 0 \;.
\]
On the other hand, if $C = \{0\}$, then
\[
K_\hT (T) = \Z \otimes_{\Z[q]} \Z((q)) \;,
\]
with $\Z$ a $\Z[q]$-module via the specialization $q \to 1$; but
by a similar argument, the resulting tensor product again vanishes. 
\end{proof}

\noindent
The functor $K^*_\hT$ extends to an equivariant cohomology theory
on the category of $\T$-CW spaces, which sends a general (large)
object $X$ to the pro-$\Z((q))$-module
\[
\{ K^*_\hT (X_i) \; | \; X_i \in {\rm finite} \; \T-{\rm CW}
\subset X \} \;,
\]
\cite{ArM} (appendix). We can thus extend the claim above: 
\bigskip

\begin{cor} For a general $\T$-CW-space $X$, the 
restriction-to-fixedpoints map
\[
K^*_\hT (X) \to K^*_\hT (X^\T)
\]
is an isomorphism of pro-objects. Moreover, if the fixedpoint
space $X^\T$ is a finite CW-space, then the pro-object on the
left is isomorphic to the {\bf constant} pro-object on the right.
\end{cor}

\noindent
The free loopspace $LX$ of a CW-space $X$ is weakly $\T$-homotopy
equivalent to a $\T$-CW-space, by a map which preserves the fixed-point
structure \cite{LMS} (\S 1.1). \bigskip

\begin{thm}
\[
M \mapsto K_\hT (LM) := K_\Tate(M)
\]
{\bf is} a cohomology theory, after all! 
\end{thm}

\begin{remark}
{\rm This seems to be what the physicists have been trying to tell us all 
along: they probably thought (as the senior author did \cite{Mo}) that the 
formal completion was a minor technical matter, not worth making any 
particular fuss about. Of course our construction is a completion of a 
much smaller (elliptic) cohomology theory, whose coefficients are modular 
forms, with the completion map corresponding to the $q$-expansion. The 
geometry underlying modularity is still \cite{Bryl} quite mysterious.}
\end{remark}

\begin{remark}
{\rm A cohomology theory defined on finite spectra extends to a cohomology 
theory on all spectra \cite{Ad}; moreover, any two extensions are equivalent,
and the equivalence is unique up to phantom maps. For the case at hand, it 
is clear that this cohomology theory is equivalent to the formal extension 
$K((q))$, where $K$ is complex $K$-theory, with $q$ a parameter in degree 
zero. Hovey and Strickland \cite{HS} have shown that an evenly graded 
spectrum does not support phantom maps, so our cohomology theory is 
uniquely equivalent, as a homotopy functor, to $K((q))$. \bigskip

\noindent
However, there is more to our construction than a simple homotopy
functor: it comes with a natural (fixed-point) orientation, which defines
a systematic theory of Thom isomorphisms. In the terminology of \cite{AHS},
it is represented by an {\bf elliptic} spectrum, associated to the 
Tate curve over $\Z((q))$; its natural orientation is defined by the 
$\sigma$-function of \S 4.} \end{remark} 
\bigskip

\begin{remark}
{\rm Let $H_{\T}$ denote $\T$-equivariant singular Borel cohomology with
rational coefficients; then $H_{\T}^*(pt) = \Q[t]$, where $t$ has degree two.
We have a localization theorem
\[ 
H^*_{\T}(X)[t^{-1}] = H^*_{\T}(X^{\T})[t^{-1}] = H^*(X^{\T})\otimes 
\Q[t^{\pm}] 
\]
for finite complexes, and can therefore play the same game as before, and 
observe that there is a unique extension $\hat{H}_{\T}( - )$ of 
$H_{\T}( - )[t^{-1}]$ to all $\T$-spectra, with a localization 
theorem valid for an arbitrary $\T$-spectrum $X$:
\[ 
\hat{H}^*_{\T}(X) = \hat{H}^*_{\T}(X^{\T}) \;. 
\] 
Jones and Petrack \cite{JP} have constructed such a theory over the real
numbers, together with the analogous fixedpoint orientation - which in
their case is (a rational version of) the $\hat A$-genus. \bigskip 

\noindent
Results of this sort (which relate oriented equivariant cohomology theories 
on free loopspaces to cohomology theories on the fixedpoints, with related 
(but distinct) formal groups, is part of an emerging understanding of what 
homotopy theorists call `chromatic redshift' phenomena, cf. \cite{AM,AMS,T}.}
\end{remark}

\begin{remark}
{\rm Our completion of equivariant $K$-theory is the natural repository 
for characters of {\bf positive-energy} representations of loop groups; 
it is {\bf not} preserved by the orientation-reversing involution $\lambda 
\mapsto \lambda^{-1}$ of $\T$. It is in some sense a {\bf chiral} 
completion.} 
\end{remark}

\begin{remark}
{\rm The completion theorem above is a specialization of Segal's 
original localization theorem \cite{Se}, which says that $K_\T(X)$, 
considered as a sheaf over the multiplicative groupscheme Spec $K_\T = 
{\mathbb G}_m$ (cf. \cite{Ro}), has for its stalks over generic (ie 
nontorsion) points, the $K$-theory of the fixed point space. The Tate point
\[
{\rm Spec} \; \Z((q)) \to {\rm Spec} \; \Z[q^\pm]
\]
is an example of such a generic point, perhaps too close to zero (or 
infinity) to have received the attention it seems to deserve. \bigskip

\noindent
In fact, we can play a similar game for any oriented theory $E_{\T}$. Let 
${\mathcal E}$ denote the union of the one-point compactifications of all 
representations of $\T$ which do not contain the trivial representation 
(cf. \cite{LMS}). Then the theory $E_\T \wedge {\mathcal E}$ satisfies 
a strong localization theorem.
} 
\end{remark}

\section{Toward Pontrjagin-Thom duality}

\noindent
One might hope for a construction which associates to a map $f: M \to N$
of manifolds (with suitable properties), a morphism
\[
Lf^{PT} : LN^{-\bT LN} \to LM^{-\bT LM}
\]
of prospectra. This seems out of reach at the moment, but some of the 
constructions sketched above can be interpreted as partial results in
this direction. \bigskip

\noindent
In particular, there is at least a {\bf cohomological} candidate for
a PT dual
\[
j^{PT} : LM^{-\bT LM} \to M^{-TM}
\]
to the fixedpoint inclusion $j: M \subset LM$. To describe it, we should
first observe that there is a morphism 
\[
M^{-TLM} \to L_!M^{-TLM}
\]
of prospectra, constructed by pulling back the tangent bundle of $LM$
along the fixedpoint inclusion. To be more precise we need to note that
the $\T$-fixedpoints of the thickened loopspace is a bundle $L_!M^\T \to M$ 
with contractible fiber; the choice of a section defines a composition
\[
\tilde j: M^{-TLM} \sim ((L_!M)^\T)^{-TLM} \to L_!M^{-TLM} \;;
\]
of course $M^{-TLM} = (M^{-TM})^{-\nu}$, where $\nu = T_M \otimes 
(\C[q^\pm]/\C)$ is the normal bundle described in \S 2. \bigskip

\noindent
Now according to the localization theorem above, the map induced on $K_\hT$
by $\tilde j$ is an isomorphism, so it makes sense to define
\[
j^! := (\tilde j^*)^{-1} \circ \phi^{-1}_\nu : K_\hT(M^{-TM}) \to 
K_\hT(M^{-TM-\nu}) \to K_\hT(LM^{-\bT LM}) \;,
\]
where $\phi_\nu$ is the Thom pro-isomorphism associated to the filtered
vector bundle $\nu$. \bigskip

\noindent
A good general theory of PT duals would provide us with a commutative 
diagram
\[
\begin{CD} LN^{-\bT LN} &@> {Lf^{PT}} >>& LM^{-\bT LM} \\@VV 
{j_N^{PT}} V &&@VV {j_M^{PT}} V\\ N^{-TN} &@> {f^{PT}} >>& M^{-TM} \;,
\end{CD}
\]
so it follows from the constructions above that 
\[
Lf^! := j^!_N \circ f^! \circ (j^!_M)^{-1}
\]
defines a formally consistent theory of PT duals for $K_\hT$. \bigskip

\noindent
For example, (the evaluation at 1 of) the composition
\[
K_\Tate(M) = K_\hT(LM) \cong K_\hT(LM^{-\bT LM}) \to K_\hT(M^{-TM}) 
\to K_\T(S^0) \cong K_\Tate(pt) 
\]
is (the $q$-expansion of) the Witten genus; more generally, our 
{\it ad hoc} construction $Lf^!$ for $K_\hT$ agrees with the covariant
construction $f^\dagger$ defined for the underlying manifold by the
fixedpoint (or $\sigma$) orientation. \bigskip

\noindent
The {\bf un}completed $K$-theory $K_\T(LM^{-{\bf T}LM})$ is also 
accessible, through the spectral sequence of a colimit, but our 
understanding of it is at an early stage. It is of course not a
cohomological functor of $M$, but it does not seem unreasonable
to hope that some of its aspects may be within reach through
similar PT-like constructions. \bigskip

\noindent
These fragmentary constructions suffice to show that $K_\hT(LM^{-\bT LM})$ 
has enough of a Frobenius (or ambialgebra) structure to define a 
two-dimensional topological field theory, which assigns to a closed 
surface of genus $g$, the class
\[
 {\pi}^{\dagger} \; \epsilon_\T(TM)^g \in \Z((q)) \;,
\]
with $\epsilon_\T$ the characteristic class defined in \S 4, and $\pi^\dagger$
the pushforward of $M$ to a point defined by the fixed-point ($\sigma$) 
orientation. When $g = 0$ this is the Witten genus of $M$, and when $g=1$ 
it is the Euler characteristic. \bigskip

\noindent
Finally: our construction is, from its beginnings onward, formulated in
terms of {\bf closed} strings. Stolz and Teichner \cite{ST} have produced a 
deeper approach to a theory of elliptic objects, which promises to 
incorporate interesting aspects of {\bf open} strings as well. However, 
their theory is in some ways quite complicated; and our hope is that 
their global theory, combined with the quite striking computational 
simplicity of the very local theory sketched here, will lead to something 
{\bf really} interesting.

\section{Appendix: The Tits Building of a loop group} \bigskip

\noindent
Let $G$ be a simply-connected compact Lie group of rank $n$. Let $LG$ denote 
the loop group of $G$. For the sake of convenience, we will work with a 
smaller (but equivalent) model for $LG$, which we now describe: \bigskip

\noindent
Let $G_{\C}$ be the complexification of the group $G$. Since $G_{\C}$ has the 
structure of a complex affine variety, we may define the group 
$L_{alg}G_{\C}$ to be the group of polynomial maps from $\C^*$ to $G_{\C}$. 
Let $L_{alg}G$ be the subgroup of $L_{alg}G_{\C}$ consisting of maps taking 
the unit circle into $G \subset G_{\C}$. The inclusion $L_{alg}G \subset LG$ 
is a homotopy equivalence in the category of $\T$-spaces. We begin by making 
our constructions with $L_{alg}G$. We then use these to draw conclusions 
about the (smooth) Tits building ${\bf A}(LG)$ \bigskip

\noindent
Fix a maximal torus $T$ of $G$, and let $\alpha_i$, $1 \leq i \leq n$ be a 
set of simple roots. We let $\alpha_0$ denote the highest root. Each root 
$\alpha_i$, $0 \leq i \leq n$ determines a compact subgroup $G_i$ of $G$. 
More explicitly, $G_i$ is the semisimple factor in the centralizer of the 
codimension one subtorus given by the kernel of $\alpha_i$. Each $G_i$ may 
be canonically identified with $SU(2)$ via an injective map 
\[ 
\varphi_i : SU(2) \longrightarrow G \;. 
\]
We use these groups $G_i$ to define corresponding compact subgroups $G_i$ 
of $L_{alg}G$ as follows:
\[ 
G_i = \{ z \mapsto \varphi_i\left(\begin{matrix} a & b \\ c & d 
\end{matrix}\right) \} \quad \quad \quad \mbox{if} \quad  \left( 
\begin{matrix} a & b \\ c & d \end{matrix}\right) \in SU(2), \quad i>0 \]
\[ G_0 = \{ z \mapsto \varphi_0\left(\begin{matrix} a & bz \\ cz^{-1} & d 
\end{matrix}\right) \}\quad \mbox{if} \quad  \left( \begin{matrix} a & b \\ 
c & d \end{matrix}\right) \in SU(2), \quad \quad \quad 
\]
\begin{remark}\label{rot}
{\rm Note that each $G_i$ is a compact subgroup of $L_{alg}G$ isomorphic to 
$SU(2)$. Moreover, $G_i$ belongs to the subgroup $G$ of constant loops if 
$i \geq 1$. The circle group $\T$ preserves each $G_i$, acting trivially 
on $G_i$ for $i \geq 1$, and nontrivially on $G_0$.}
\end{remark}

\begin{defn} 
{\it For any proper subset $I \subset \{0,1,\ldots,n\}$, define the 
parabolic subgroup $H_I \subset L_{alg}G$ to be the group generated by $T$ 
and the groups $G_i$, $i \in I$. For the empty set, we define $H_I$ to be 
$T$. It follows from \ref{rot} that each $H_I$ is preserved under the 
action of $\T$.} 
\end{defn}

\begin{remark}
{\rm The groups $H_I$ are compact Lie \cite{Mi}. Moreover, $H_I$ is 
isomorphic to its image in $G$, under the evaluation map $ev(1) : L_{alg}G 
\rightarrow G$. Notice that for $I = \{1,\ldots,n\}$, $H_I = G$. Notice 
also that $\T$ acts nontrivially on $H_I$ if and only if $0 \in I$.}  
\end{remark}

\noindent
We are now ready to define the Tits building ${\bf A}(L_{alg}G)$. 
\begin{defn}
{\it Let ${\bf A}(L_{alg}G)$ be the homotopy colimit:
\[ 
{\bf A}(L_{alg}G) = {\rm hocolim}_{I \in \mathcal{C}} \; L_{alg}G/H_I 
\]
where $\mathcal{C}$ denotes the poset category of proper subsets of 
$\{0,1,\ldots,n\}$.} 
\end{defn}

\noindent
We now come to the main theorem: 

\begin{thm}\label{contr}
The space ${\bf A}(L_{alg}G)$ is $\T \tilde{\times} L_{alg}G$-equivariantly 
contractible. In other words, given a compact subgroup $K \subset \T 
\tilde{\times} L_{alg}G$, then the fixed point space ${\bf A}(L_{alg}G)^K$ 
is contractible.
\end{thm}

\begin{proof}
A proof of the contractibility of ${\bf A}(L_{alg}G)$ was given in 
\cite{Mi}. We use some of the ideas from that paper, but our proof is 
different in flavour. \bigskip

\noindent
Mitchell expresses the space ${\bf A}(L_{alg}G)$ as the following:
\[
{\bf A}(L_{alg}G) = (L_{alg}G/T \times \Delta)/\sim
\] 
where $\Delta$ is the $n$-simplex, and $(aT,x) \sim (bT,y)$ if and only if 
$x=y \in \stackrel{\circ}{\Delta}_I$ and $aH_I= bH_I$. Here we have indexed 
the walls of $\Delta$ by the category $\mathcal{C}$, and denoted the 
interior of ${\Delta}_I$ by $\stackrel{\circ}{\Delta}_I$. \bigskip

\noindent
Let $L_{alg}\gG \oplus \R d$ be the Lie algebra of the extended loop group 
$\T \tilde{\times} L_{alg}G$. Consider the affine subspace 
\[ 
{\bf A} = L_{alg}\gG + d \subset L_{alg}\gG \oplus \R d \;.
\]
The adjoint action of $L_{alg}G$ on ${\bf A}$ is given by
\[ 
Ad_{f(z)}(\lambda(z) + d) = Ad_{f(z)}(\lambda(z)) + zf'(z)f(z)^{-1} + d 
\]
This action extends to an affine action of $\T \tilde{\times} L_{alg}G$. The 
identification of ${\bf A}$ with ${\bf A}(L_{alg}G)$ is given as follows. 
Let $\Delta$ be identified with the affine alcove:
\[ 
\Delta = \{ (h+d) \in Lie(T) + d \, \, | \quad \alpha_i(h) \geq 0, \, i>0, 
\, \alpha_0(h) \leq 1 \} \;.
\]
General facts about Loop groups \cite{Ka, Mi} show that the surjective map
\[ 
L_{alg}G \times \Delta \longrightarrow A, \quad (f(z),y) \mapsto 
Ad_{f(z)}(y) 
\]
has isotropy $H_I$ on the subspace $\Delta_I$. Hence it factors through a 
$\T \tilde{\times} L_{alg}G$-equivariant homeomorphism between ${\bf A}
(L_{alg}G)$ and the affine space ${\bf A}$. Notice that any compact subgroup 
$K \subset \T \tilde{\times} L_{alg}G$ admits a fixed point on ${\bf A}
(L_{alg}G)$. Hence, the space ${\bf A}(L_{alg}G)^K$ is also affine. This 
completes the proof. 
\end{proof}

\noindent
We now define the smooth Tits building

\begin{defn}
{\it Let ${\bf A}(LG)$ be the homotopy colimit:}
\[ 
{\bf A}(LG) = {\rm hocolim}_{I \in \mathcal{C}} \; LG/H_I = LG 
\times_{L_{alg}G} {\bf A}(L_{alg}G) \;.
\]
\end{defn}

\noindent
It is clear from the proof of the above theorem, that ${\bf A}(LG)$ is $\T 
\tilde{\times} LG$-equivariantly homeomorphic to the affine space $LG 
\times_{L_{alg}G}{\bf A}$ which is homeomorphic to the affine space 
${\mathcal A}(S^1 \times G)$ of connections on the trivial bundle $S^1 
\times G$. This shows that:
 
\begin{thm} 
The smooth Tits building ${\bf A}(LG)$ is $\T \tilde{\times} 
LG$-equivariantly contractible. 
\end{thm}

\noindent
Recall \cite{Mi} that ${\bf A}(L_{alg}G)$ is homeomorphic to the space
\[  
{\mathcal S}_{alg} = \{ g(t):[0,1] \rightarrow G \quad |\, g(t) = 
f(e^{2\pi it})\cdot exp(tX);\, f(z) \in \Omega_{alg}G, X \in Lie(G) \} \;
\]
We have a corresponding smooth version
\[ 
{\mathcal S} = \{ g(t):[0,1] \rightarrow G  \quad |\, g(t) = 
f(e^{2\pi it})\cdot exp(tX);\, f(z) \in \Omega G \} = LG \times_{L_{alg}G}
{\mathcal S}_{alg}  \;
\]
which is clearly homeomorphic to ${\bf A}(LG)$. It remains to identify 
${\mathcal S}$ with the space 
\[ 
{\mathcal S} = \{ g(t):\R \rightarrow G, \, \, g(0)=1, \, \, g(t+1) = g(t)
\cdot g(1) \} 
\;;
\]
This is straightforward, and is left to the reader. 

\newpage

\bibliographystyle{amsplain}

\end{document}